\documentclass[preprint]{elsarticle}
\pdfoutput=1
\usepackage{amsmath}
\usepackage{amssymb}

\renewcommand{\vec}[1]{{\mathbf #1}}
\newcommand{\R}{{\mathbb R}}

\newcommand{\Q}{{\mathbb Q}}

\newcommand{\eps}{{\varepsilon}}

\newtheorem{lemma}{Lemma}

\newtheorem{theorem}[lemma]{Theorem}
\newtheorem{corollary}[lemma]{Corollary}

\begin{document}

\begin{frontmatter}

\title{Reconstructions in Ultrasound Modulated Optical Tomography}
\author{Moritz Allmaras}
\ead{allmaras@math.tamu.edu}
\author{Wolfgang Bangerth\corref{cor1}}
\ead{bangerth@math.tamu.edu}
\address{Department of Mathematics, Texas A\&M University,\\ College
Station, TX 77843, USA}
\cortext[cor1]{Corresponding author. Address: Department of Mathematics, Texas
A\&M University, 3368 TAMU, College Station, TX 77843, USA; Tel.: +1 979 845
6393, Fax: +1 979 862 4190}

\begin{abstract}
We introduce a mathematical model for ultrasound modulated optical tomography
and present a simple reconstruction scheme for recovering the spatially
varying optical absorption coefficient from scanning measurements with narrowly
focused ultrasound signals. Computational results for this model show that the
reconstruction of sharp features of the absorption coefficient is possible. A
formal linearization of the model leads to an equation with a Fredholm operator,
which explains the stability observed in our numerical experiments.
\end{abstract}


\end{frontmatter}

\section{Introduction}

During the last two decades, optical tomography (OT) has received significant
attention as a biomedical imaging modality. This can be attributed, in
particular, to the fact that light at optical frequencies is harmless to the
human organism and that optical properties of tissues reveal important
biological information such as angiogenesis and hypermetabolism, both of which
are well-known indicators of cancer \cite{Wang07}.  Unfortunately,
reconstruction in OT is also known to be severely ill-posed, and consequently
the sharp imaging of optical properties is all but impossible. Various
attempts to address this problem have been made. In this paper we are
interested in a hybrid imaging method called Ultrasound Modulated Optical
Tomography (UOT, \cite{Wang07}) that combines the OT procedure with
simultaneous modulation by a narrowly focused ultrasound beam in order to
alleviate the instability of OT reconstructions.  The idea is to combine the
good tumor specificity of OT with the high spatial resolution of ultrasound
imaging. This approach utilizes the experimentally observed interaction
between ultrasound and light propagation in tissue \cite{LM95,Wang07}. In UOT,
a coherent light source irradiates the tissue sample and causes interference
patterns to form on the surface of the object, so-called speckles.  A narrowly
focused ultrasound wave is simultaneously induced in the tissue, influencing
its optical properties and thus modulating the speckle pattern with ultrasound
frequency. By measuring properties of this modulation, information about the
incident light intensity at the focus location of the ultrasound beam can be
obtained. Hence, by scanning the focus of the ultrasound wave throughout the
sample, a quantity related to the light intensity in the object's
\textit{interior} can be determined. This type of internal information is
usually not available from OT measurements due to multiple scattering of
photons in optically dense media, although there are other variants of optical
tomography that also strive to recover this information
(e.g. \cite{SKGHR03}). It can be expected that this additional knowledge can
help in stabilizing the inversion process and render it substantially less
ill-posed than the original OT problem. For the UOT model we present in this
paper, numerical experiments and an initial analysis suggest that this
intuition is justified.

The literature contains a number of models that address the UOT
technique, see for example \cite{LM95,KLZG97,Wan01,Nam02,Wan03,Wang07}.
Most of them describe the coupling between ultrasound and light in terms of
stochastic quantities, which permits particle-based simulations of the
light intensity modulation effect caused by the ultrasound wave. On the other
hand, for
optical imaging in turbid media at a depth of several centimeters, photon
intensities can be accurately modeled by the diffusion limit. Under
certain assumptions, this allows us to formulate a model
for the UOT procedure based on a parameter identification problem for a set of
coupled diffusion-type partial differential equations. This model, along with a
description of the measurements is presented in Section~\ref{sec:model}.
In Section~\ref{sec:algorithms}, we outline a simple algorithm that can be used
to reconstruct the spatially varying absorption coefficient from UOT
measurements with focused ultrasound signals. Examples of the resulting reconstructions for numerical
phantoms are provided in Section~\ref{sec:implementation}. In Section
\ref{sec:stability}, we formally linearize our model and obtain an equation
that relates perturbations in the absorption coefficient to those in the
measurements
by a Fredholm operator acting between appropriate Sobolev spaces. This
provides a partial explanation to the stable reconstruction observed in
our numerical experiments. The last section contains final remarks and
conclusions.

\section{Mathematical model}\label{sec:model}

A detailed description of the physical underpinnings of the UOT procedure can
be found, for instance, in \cite[Ch. 13]{Wang07}. We give a brief description of
the set-up here.

Let the object of interest occupy the domain $\Omega \subset \R^3$. The
internal optical properties in the diffusion limit are described by the reduced
scattering coefficient $\mu_s'$ and the absorption coefficient $\mu_a$. For
imaging soft tissues, it is common to assume $\mu_s'$ roughly equal to a known
constant throughout $\Omega$, while the spatially varying absorption
$\mu_a(x)$, $x\in\Omega,$ represents the target of reconstruction. It is also
assumed that the
tissue of interest is turbid (highly scattering), so that $\mu_a(x)\ll
\mu_s'$. It is known that in such media, the light intensity $u(x)$ inside
$\Omega$ can be accurately described by the diffusion approximation (e.g., \cite{Cha60,SKGHR03}).

It has been shown experimentally that coherent light can be modulated by an ultrasound
field inside the turbid medium \cite{LM95}. Various explanations have been put
forward for this effect \cite{Wan01}.

The experimental setup in UOT involves
dealing with the time dependent light intensity of individual speckles.
The model presented below is derived under two assumptions, which are satisfied
in standard UOT applications \cite{Wang07}:
 \begin{itemize}
   \item Weak scattering assumption: The optical wavelength is much shorter
than the mean free path.
   \item Weak ultrasound modulation assumption: The ultrasound-induced change
in the optical path length is much less than the optical wavelength.
 \end{itemize}

The measured signal is the autocorrelation function \cite{LM95} at a detector
location $\eta\in\partial\Omega$
\[
  G_1(\eta,\tau) = \left<E(\eta,t+\tau)E^\ast(\eta,t)\right>_t,
\]
where angle brackets denote averaging over time, and the electric field $E$ is
related to the light intensity $I$ as $I(\eta,t)=|E(\eta,t)|^2$.  It has been
shown experimentally \cite{LM95} that over time scales $\tau\gg 1\mu s$
coherence of the exiting light is lost, i.e. $G_1(\eta,\tau)\rightarrow 0$ as
$\tau\rightarrow\infty$, due to the Brownian motion of scatterers. However, on
short time scales on the order of the period of the ultrasound field --
i.e. the regime we are interested in --, $G_1(\eta,\tau)$ has been observed to
oscillate at the ultrasound frequency.  We will therefore neglect contributions from the
Brownian motion of scatterers since it is unrelated to the ultrasound field. In the
absence of an ultrasound field, and on these time scales, we would then have
$G_1(\eta,\tau)=\textrm{const}$. In the following, we will derive expressions
for $G_1$ and, in particular, its \textit{modulation depth}, i.e. the
magnitude of the oscillation of $G_1$ at the ultrasound frequency. We will
then relate these quantities to solutions of partial differential equations
that we will use for our reconstruction scheme.

\paragraph{A path integral model}
For a point source of unit strength at a location $\sigma$, and a detector
measuring photons exiting the domain at $\eta\in\partial\Omega$, we can
write
\begin{equation*}
  G_1(\sigma,\eta,\tau)
  = P^\partial \bar G(\sigma,\eta,\tau),
  \qquad
  \bar G(\sigma,\eta,\tau)
  = \sum_{s=s(\sigma,\eta)} P_s \left<E_s(t+\tau)E^\ast_s(t)\right>_t
\end{equation*}
where the sum extends over all paths $s$ that connect source $\sigma$ and
boundary location $\eta$. $P_s$ is the fraction of the incident
intensity that scatters along $s$ multiplied by the probability of a
photon not getting absorbed along this path. $P^\partial$ is the probability
that a photon that makes it to a
point $\eta$ on the boundary is able to cross the boundary from tissue into
the detector. $E_s$ then denotes the
phase of the electric field at $\eta$ of photons following path
$s$. Consequently, $\left<E_s(t)E^\ast_s(t)\right>_t=1$.

Consider now the situation in which the ultrasound field $p(x,t)$ induces
phase shifts $d\phi(x,t)$ on all paths along an infinitesimal path element
$ds(x)$. As shown in \cite{Wan01}, such phase shifts can be induced both by
the periodic motion of scatterers in the ultrasound field as well as by the
modulation of the index of refraction by the pressure field. We then have
\begin{eqnarray*}
  \left<E_s(t)E^\ast_s(t+\tau)\right>_t
  &=&
  \left< \exp\left( -i \int_s
      \frac{d\phi(x,t)}{ds} \; ds \right)\right>_t
  \\
  &\approx&
  \exp\left( -\frac 12
    \left<
    \left[ \int_s
      \frac{d\phi(x,t)}{ds} \; ds
    \right]^2
    \right>_t \right),
\end{eqnarray*}
where integrals are assumed to be along a path $s$ from $\sigma$ to $\eta$.
By computing how the index of refraction and the phase shifts induced by
scatterer movement depend on an ultrasound pressure field with frequency
$\omega_a$, we can use the results in \cite{Wan01} to write above expression
as
\begin{eqnarray*}
  \left<E_s(t)E^\ast_s(t+\tau)\right>_t
  &=
  \frac 1{|s|}
  \int_s \exp\left[
    -\alpha |p(x)|^2 (1-\cos \omega_a\tau)
  \right]
  \; ds,
\end{eqnarray*}
where $\alpha$ is a proportionality
constant and $|s|$ is the length of path $s$. (Note in particular that the
proportionality to the square of the
pressure has also been observed experimentally, see \cite{LM95}.)
Consequently,
\begin{equation*}
  G_1(\sigma,\eta,\tau)
  =
  P^\partial
  \sum_{s=s(\sigma,\eta)}
  P_s
  \frac 1{|s|}
  \int_s \exp\left[
    -\alpha |p(x)|^2 (1-\cos \omega_a\tau)
  \right]
  \; ds.
\end{equation*}
As has been shown experimentally \cite{LM95}, the temporal variation of the
exponent is relatively small. We can therefore approximate
\begin{eqnarray}
  G_1(\sigma,\eta,\tau)
  &=
  P^\partial
  \sum_{s=s(\sigma,\eta)}
  P_s
  \left[
    1
    -
    \frac {\alpha}{|s|}\int_s
    |p(x)|^2 (1-\cos \omega_a\tau)
    \; ds
  \right].
\end{eqnarray}
It follows that we can write the autocorrelation function as the sum of two
terms:
\begin{eqnarray*}
  G_1(\sigma,\eta,\tau)
  &=
  G_1(\sigma,\eta,0)
  -
  \alpha
  P^\partial
  \sum_{s=s(\sigma,\eta)}
  P_s
  \frac 1{|s|}
  \int_s
  |p(x)|^2 (1-\cos \omega_a\tau)
  \; ds.
\end{eqnarray*}
The first of these is the time average light intensity, whereas the second is
the temporal variation of the autocorrelation function due to the ultrasound
field. To first order in the small parameter $\alpha$, this expression equals
\begin{eqnarray*}
  G_1(\sigma,\eta,\tau)
  &=
  G_1(\sigma,\eta,0)
  -
  \alpha
  P^\partial
  \int_\Omega
  \bar G(\sigma, x, 0)
  |p(x)|^2
  \bar G(x, \eta, 0)
  \; dx
  \; (1-\cos \omega_a\tau).
\end{eqnarray*}

Finally, if light is incident with an intensity $S(\sigma)$ at source
positions $\sigma\in\partial\Omega$, the overall autocorrelation function at
detector location $\eta$ can be written as
\begin{eqnarray}
  G_1(\eta,\tau)
  &=
  \int_{\partial\Omega} S(\sigma) G_1(\sigma,\eta,\tau) \; d\sigma.
\end{eqnarray}
Using the previous equation, and defining the time averaged light intensity
$u(x)=\int_{\partial\Omega} S(\sigma) \bar G(\sigma, x,0) \; d\sigma$
for all $x\in \Omega\cup\partial\Omega$, we can then write
\begin{eqnarray}
  G_1(\eta,\tau)
  &=&
  P^\partial
  u(\eta)
  -
  \alpha
  P^\partial
  \int_\Omega
  u(x) |p(x)|^2
  \bar G(x, \eta, 0)
  \; dx
  \; (1-\cos \omega_a\tau)
  \\
  &=&
  P^\partial
  \left[
    u(\eta)
    -
    v(\eta) (1-\cos \omega_a\tau) \right],
\end{eqnarray}
where
\begin{equation}
  \label{eq:def-v}
  v(\eta)=\alpha \int_\Omega u(x) |p(x)|^2 \bar G(x, \eta, 0) \;
  dx.
\end{equation}
This representation of the correlation as a sum of a time averaged photon flux
plus a temporally variable term has given rise to the name \textit{tagged
  photons} to denote $v(x)$. Equation \eqref{eq:def-v} makes it clear
that tagged photons originate at the site $x$ of interaction of the
steady-state light field $u(x)$ and ultrasound field $p(x)$. However,
since $G_1(\eta,\tau)$ is not a photon flux but a correlation function, we
will not use this term any further.

In our reconstruction algorithm below, we will assume that the amplitude $
v(\eta)$ of the temporal variation of $G_1(\eta,\tau)$ -- i.e. the
\textit{modulation depth} -- is the measured signal. While the time average
$u(\eta)$ can also be measured, using it for inversion leads to the
diffuse optical tomography problem that is known to be severely ill-posed.

\paragraph{A partial differential equation model}
For our reconstruction algorithms, we would like to relate our signal
$v(\eta)$ to the solution of a partial differential equation. To this end,
note
that $\bar G(x,y,0) = \sum_{s=s(x,y)} P_s$ is the time average
probability that a photon starting at
$x$ is found at $y$.
For the turbid medium that we
consider in this contribution, light propagation can be accurately described by
the diffusion approximation in which photons perform a random walk. The time
averaged light intensity $u(\eta)=\int_{\partial\Omega} S(\sigma)
\bar G(\sigma,\eta,0) \; d\sigma$ must then satisfy the following
equation:
\begin{equation}  \label{eq:diffusion-u}
  -\nabla\cdot D\nabla u(x)+\mu_a(x) u(x)
  =0 \qquad \textrm{in}\ \Omega,
\end{equation}
where
\begin{equation}\label{eq:D}
D=D(x)=\frac{1}{3(\mu_a(x)+\mu_s'(x))}
\end{equation}
is the diffusion coefficient. Due to the assumptions stated at the beginning
of this section,
$D\approx \frac{1}{3\mu_s'}\approx\textrm{const}$. To simplify the notation, we
set $\mu:=\mu_a$ in the rest of the text.
Equation~\eqref{eq:diffusion-u} needs to be completed by boundary conditions.
For tissue in contact with a surrounding medium, Robin-type boundary conditions
are typically chosen~\cite{GHRS02}:
\begin{equation}
  \label{eq:diffusion-u-boundary}
  2D\frac{\partial u(x)}{\partial{n}}+\gamma u(x)
  = S(x)\qquad\qquad \textrm{on}\ \partial\Omega.
\end{equation}
Here $n$ denotes the outward normal to the surface $\partial\Omega$ and
$\gamma > 0$ is a constant describing the optical refractive index mismatch
at the boundary, and is related to $P^\partial$. In particular,
the assumptions underlying the diffusion approximation imply that
$P^\partial\le \frac 12$ and $0\le\gamma \le 1$, and $\gamma=1$ if
$P^\partial=\frac 12$.

On the other hand, to represent $v(\eta)$ as the solution of a partial
differential equation, we have to consider the equation that $\bar G$
satisfies. $\bar G(x,y,0)$ is the probability that a photon originating at $x$
reaches $y$, absent an ultrasound field. For random walk models, it is
known that $\bar G(x,y,0)$ satisfies a diffusion
equation \cite{SaUe65,GaMcK72}, which in our case is
\begin{equation*}
  -\nabla\cdot D\nabla \bar G(x,y,0)+\mu(x) \bar G(x,y,0)
  = \delta(x-y) \qquad \textrm{ in }\ \Omega.
\end{equation*}
The question of boundary conditions is less clear. It is well known that if
every particle that reaches the boundary leaves the domain
(i.e. $P^\partial=1$), then the correct boundary condition to choose is
$\bar G|_{\partial\Omega}=0$. On the other hand, if all photons are reflected
and
none can leave (i.e. $P^\partial=0$), then $n \cdot \nabla
\bar G|_{\partial\Omega}=0$ is the correct boundary condition. In either of
these
two cases, $n \cdot \nabla \bar G|_{\partial\Omega}$ is the flux of particles
across the boundary.
However, we have been unable to find literature on the case
$0<P^\partial<\frac{1}{2}$ (see, however, \cite{DGCa76} for the case where
each particle that reaches the boundary is replaced by more than one new
particle, a situation that formally corresponds to the situation where the
fraction of particles that can leave the domain satisfies
$P^\partial<0$). Since intuitively, $\bar G$ denotes a photon flux,
we conjecture by way of analogy that $\bar G$ also satisfies
Robin boundary conditions
\begin{equation}
  2D\frac{\partial \bar G(x,y,0)}{\partial{n}}+\gamma \bar G(x,y,0)
  = 0 \qquad \textrm{ on }\ \partial\Omega.
\end{equation}

Under this assumption, we have that the amplitude
$v(\eta)$
(up to the constant factor $P^\partial$) of the time variation of
the autocorrelation function $G_1(\eta,\tau)$ satisfies the following boundary
value problem:
\begin{equation}\label{eq:diffusion-v}
\left\{
\begin{array}{ll}
  -\nabla\cdot D\nabla v(x)+\mu(x) v(x)
  =\alpha |p(x)|^2 u(x) & \textrm{ in }\ \Omega,
  \\
  2D\frac{\partial v(x)}{\partial{n}}+\gamma v(x)
  = 0 & \textrm{ on }\ \partial\Omega.
\end{array}
\right.
\end{equation}
Note that if we were to view $v$ as a fluence of \textit{virtual} or
\textit{tagged} photons, then our conjecture implies that the equation for
this virtual fluence has the same boundary conditions as that for the incident
fluence $u$.

\paragraph{Measurements}
In principle, the interferometric detectors for the modulation $P^\partial
v(\eta)$ visible beyond the boundary could be placed along the entire
boundary. In practice, however, we will only be able to measure at a small
number of locations. To simplify the
discussion, we will assume in the following that only a single detector
is used. More elaborate experimental setups could use multiple detectors to
suppress the effects of noise on the reconstruction.

\paragraph{The inverse problem}
We can now formulate the inverse problem addressed in this work:
{\em Assuming that for a given point $\eta\in\partial\Omega$ and a number of
  ultrasound fields $p^\xi(x)$ indexed by $\xi$, the values
\begin{equation}\label{eq:measurement}
h(\xi):=v^\xi (\eta)
\end{equation}
are known in the coupled system of equations}
\begin{equation}\label{eq:inverse_pr}
\left\{
\begin{array}{ll}
    -\nabla\cdot D\nabla u(x) +\mu(x) u(x)=0
    & \qquad\qquad \textrm{ in }\ \Omega,
    \\
    2D\frac{\partial u(x)}{\partial{n}}+\gamma u(x)
    = S(x) & \qquad\qquad \textrm{ on }\ \partial\Omega,
    \\
    -\nabla\cdot D\nabla v^\xi (x) +\mu(x) v^\xi (x)
    =\alpha |p^\xi(x)|^2 u(x) & \qquad\qquad \textrm{ in }\ \Omega,
    \\
    2D\frac{\partial v^\xi (x)}{\partial{n}}+\gamma v^\xi (x)
    = 0 & \qquad\qquad \textrm{ on }\ \partial\Omega.

                                       \end{array}
                                     \right.
\end{equation}
{\em Then recover the absorption coefficient $\mu$ inside a region of interest
  $U \subset \Omega$ with $\bar U \subset \Omega$.}

We remark that in the applications of ultrasound modulated optical tomography
available in the literature, the ultrasound pressure field $p(x)$ is always a
beam focused on a single point. In particular, the algorithm we show below
is based on the assumption of perfectly focused beams
$|p^\xi(x)|^2=\delta(x-\xi)$, although we will test
in Section~\ref{sec:elongated} how the algorithm performs on data for which
this assumption is not satisfied. The formulation above is more general in
that it allows arbitrary fields $p(x)$. An application of this includes
ultrasound pressure fields that are focused not on points but on spherical
surfaces for synthetic focusing, as mentioned in Section~\ref{sec:synthetic}.

\section{Reconstruction algorithm}
\label{sec:algorithms}

In this section, we introduce a simple algorithm that can be used to compute
numerical reconstructions for the above inverse problem. In the following, we
will make the assumption that the pressure field is perfectly focused on a
location $\xi\in\Omega$, i.e. $|p^\xi(x)|^2=\delta(x-\xi)$. As discussed
in Section~\ref{sec:elongated}, this is of course not practically feasible, so
our assumption is understood to mean that the real pressure field approximates
a perfectly focused one.

Let $G(x,y)$ be the Green's function for the diffusion
model \eqref{eq:diffusion-u}, i.e. the solution of
\begin{equation}\label{eq:green}
\left\{
\begin{array}{ll}
  -\nabla_x\cdot D\nabla_x G(x,y) +\mu(
  x) G(x,y)
  =\delta(x-y) & x\in \Omega,
  \\
  2D\frac{\partial G(x,y)}{\partial{n}}+\gamma G(x,y)
  = 0 & x\in \partial\Omega.
\end{array}
\right.
\end{equation}
Then, (\ref{eq:inverse_pr}) implies
\[
  v^\xi (x)= \alpha G(x,\vec \xi) u(\vec \xi),
\]
and thus,
\[
  h(\xi)=\alpha G(\eta, \xi)u(\xi), \qquad u(\xi)=\frac{h(\xi)}{\alpha G(\eta,
\xi)}.
\]
Substituting this expression for $u$ into the first equation of
(\ref{eq:inverse_pr}), we obtain an equation for recovering $\mu$:
\begin{equation}
  \label{eq:mu-recon}
  \mu(\xi) = \frac{[\nabla_\xi\cdot D\nabla_\xi]\left(h(\xi)/G(\eta,
\xi)\right)}
             {h(\xi)/G(\eta, \xi)}.
\end{equation}
The apparent difficulty in using this formula for reconstruction is that it is
implicit in $\mu$ since both $D$ and the Green's function $G$ depend on the
absorption. However, we can construct the following natural iterative
scheme for \eqref{eq:mu-recon}:
\begin{itemize}
  \item \textbf{Initial step}: Using an initial guess $\mu^0$ for the absorption
coefficient (e.g. $\mu^0=\textit{const}$), compute the corresponding
Green's function numerically, and apply formula \eqref{eq:mu-recon} to find a
new approximation $\mu^1$ for the absorption.
\item \textbf{Iterative step:} Using the current approximation $\mu^k$,
re-compute Green's function and $D$ and apply formula \eqref{eq:mu-recon} to
find an updated absorption coefficient $\mu^{k+1}$.
\end{itemize}
We do not consider the convergence properties of this scheme here, but note
that in our numerical tests presented below the iterates converged reliably,
albeit not very rapidly.

\section{Numerical implementation}
\label{sec:implementation}

Implementation of the algorithm outlined above requires the following steps:
\begin{itemize}
\item Simulation of the forward model to generate synthetic measurements,
\item repeated computation of the Green's function $G(x,y)$ for
  equation (\ref{eq:green}),
\item repeated evaluation of the iteration formula (\ref{eq:mu-recon}).
\end{itemize}
These steps are discussed in the following
subsections. In this work, we only consider measurements obtained by forward
calculations from mathematical phantoms, rather than actual experimental data.
All computations were done in $2$D, although they can be readily carried over to
$3$D. For the finite element calculations involved in
the reconstruction scheme, the Open Source finite element library deal.II
\cite{BHK07,BK99m} was used.

\subsection{Forward simulations}
In order to generate the measurements $h(\xi)$ (see (\ref{eq:measurement})),
we need to compute the solution $u(x), v^\xi(x)$ of the forward
problem (\ref{eq:inverse_pr}) for a set of given data $D, \mu, S$ (diffusion
coefficient, absorption coefficient, incoming light flux) and an ultrasound
signal focused at the point $\xi\in U$. Then, evaluating $v^\xi$ at the detector
location $\eta$, we obtain the measurement value $h(\xi)$.

\subsubsection{Computational setting.}\label{sec:setting}
\begin{figure}[htbp]
  \centering
  \includegraphics[scale=0.7]{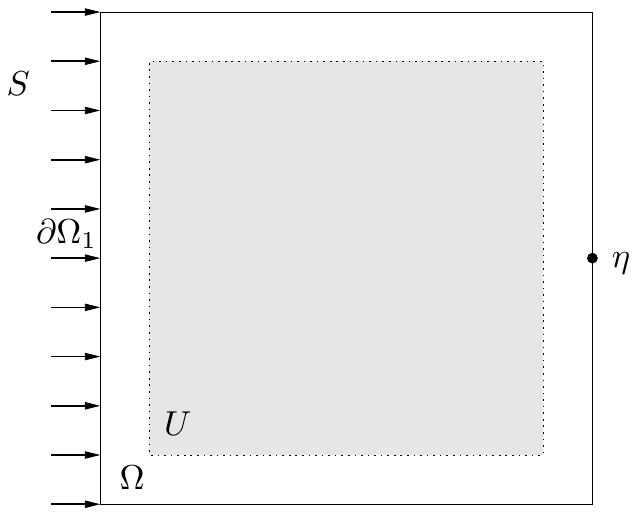}
  \caption{Setting of numerical experiments: Domain $\Omega$, area of interest
    $U$, incident light source $S(x)$ on the left, and
    detector point $\eta$ on the right.}
  \label{fig:setting}
\end{figure}

We take $\Omega$ to be the square $[0,5\rm{cm}]^2$, which approximately
corresponds to the relevant dimensions in practical applications.
For the boundary light source $S$ in \eqref{eq:diffusion-u-boundary},
$\partial\Omega$ is split into $\partial\Omega_1 = \{x\in \partial\Omega:
x_1=0\}$ and $\partial\Omega_2 = \partial\Omega \setminus \partial\Omega_1$.
Constant illumination is assumed on $\partial\Omega_1$ and no photons are
injected on $\partial\Omega_2$:
\begin{equation}
S(x)=\left\{
\begin{array}[2]{l@{\quad}l}
   1 & \textrm{for}\quad x \in \partial\Omega_1,\\
   0 & \textrm{for}\quad x \in \partial\Omega_2.
\end{array}\right.
\end{equation}
The modulation depth is measured at a single detector
location $\eta=(5\rm{cm}, 2.5\rm{cm})$. This layout is depicted in
Fig.~\ref{fig:setting}.

\subsubsection{Incident light field.}\label{sec:incident}
Since in our model the incident light intensity $u$ is independent
of the shape and location of the ultrasound waves in the tissue, $u$ only needs
to be computed once. For this computation, a finite element approximation to $u$
is constructed on a regular rectangular grid using $\Q_1$ finite
elements \cite{BS02}, solving equations
\eqref{eq:diffusion-u}--\eqref{eq:diffusion-u-boundary}. The left panel of
Fig.~\ref{fig:u-orig-v} shows $u$ for the case of a constant absorption
coefficient $\mu$.

\subsubsection{Ultrasound field.}\label{sec:us-functions}
In our numerical examples, we use Gaussian-shaped synthetic ultrasound signals:
\begin{equation}\label{eq:field}
 p(x) = C\exp\bigg(-\sum_{j=1}^d
\frac{|x_j|^2}{\sigma_j^2}\bigg),
\end{equation}
where $C$ is a normalization constant.
By choosing different variances $\sigma_j^2$ we can model varying focusing
properties of such pressure field.

To simulate scanning of the ultrasound focus, focusing points $\{\vec
\xi^i,\; i=1,\ldots,N\}$ are placed at the vertices of a square grid covering
the area of interest, here chosen as the square $U =
[0.5\rm{cm},4.5\rm{cm}]^2\subset\Omega$. For each $i$ we
then construct a signal $p^{\xi_i}(x)$ focused at $\xi_i$ by setting
\[
p^{\xi_i}(x) := p(x-\xi_i).
\]
To simplify notation we set $v^i:=v^{\xi^i}$ and $p^i:=p^{\xi^i}$.

\subsubsection{Modulated light field and measurements.}\label{sec:modulated}
\begin{figure}[htbp]
  \centering
  \begin{minipage}[b]{0.32\textwidth}
    \includegraphics[width=\textwidth]{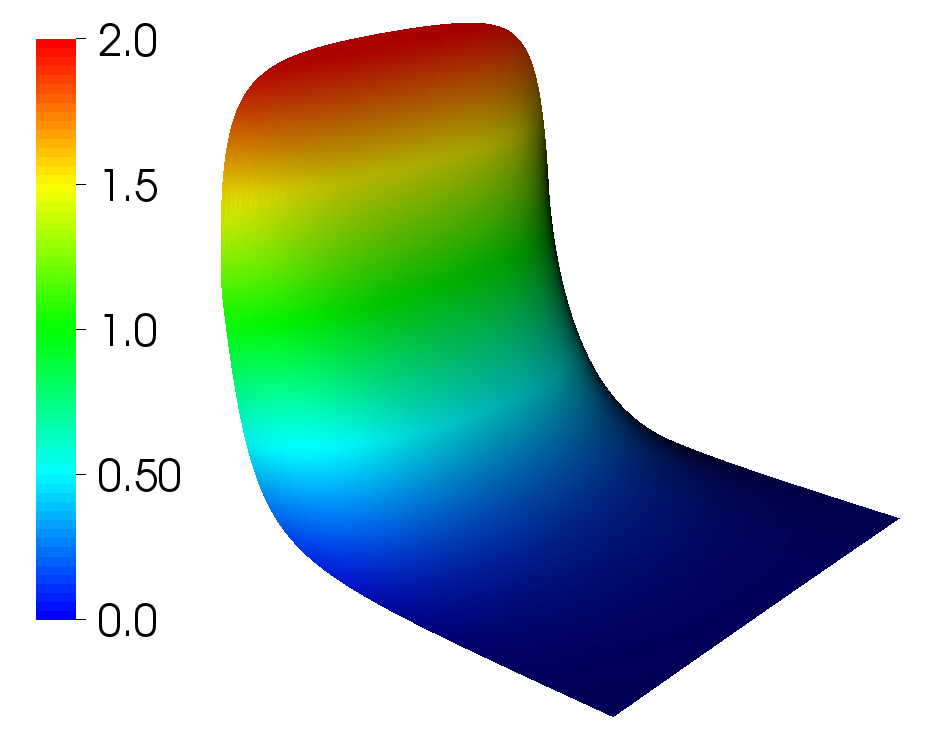}
  \end{minipage}
  \begin{minipage}[b]{0.32\textwidth}
    \includegraphics[width=\textwidth]{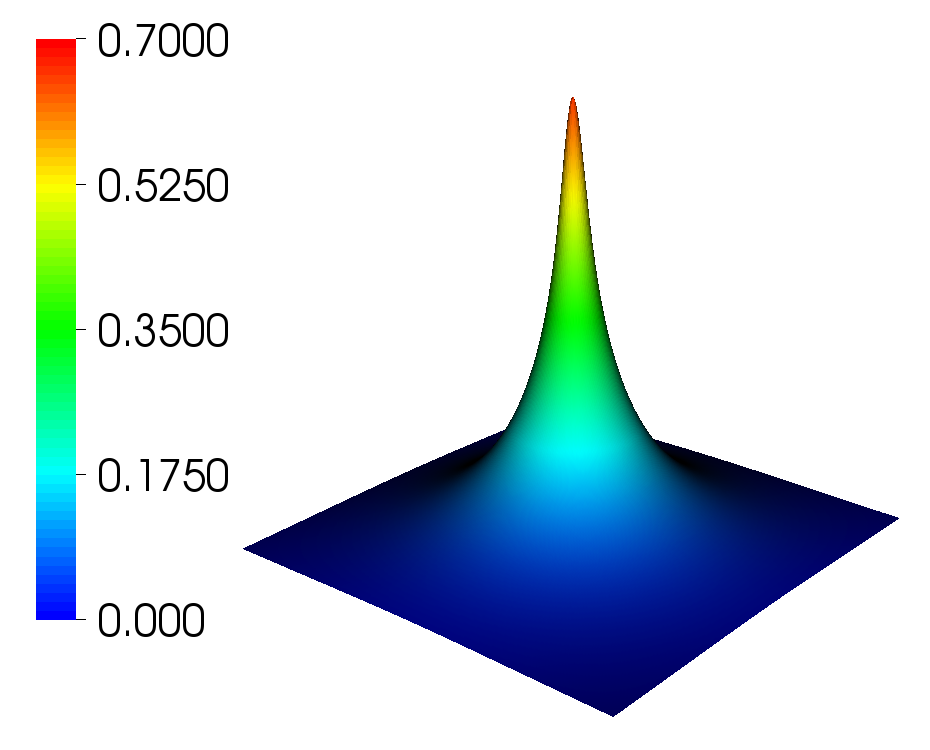}
  \end{minipage}
  \begin{minipage}[b]{0.32\textwidth}
  \includegraphics[width=\textwidth]{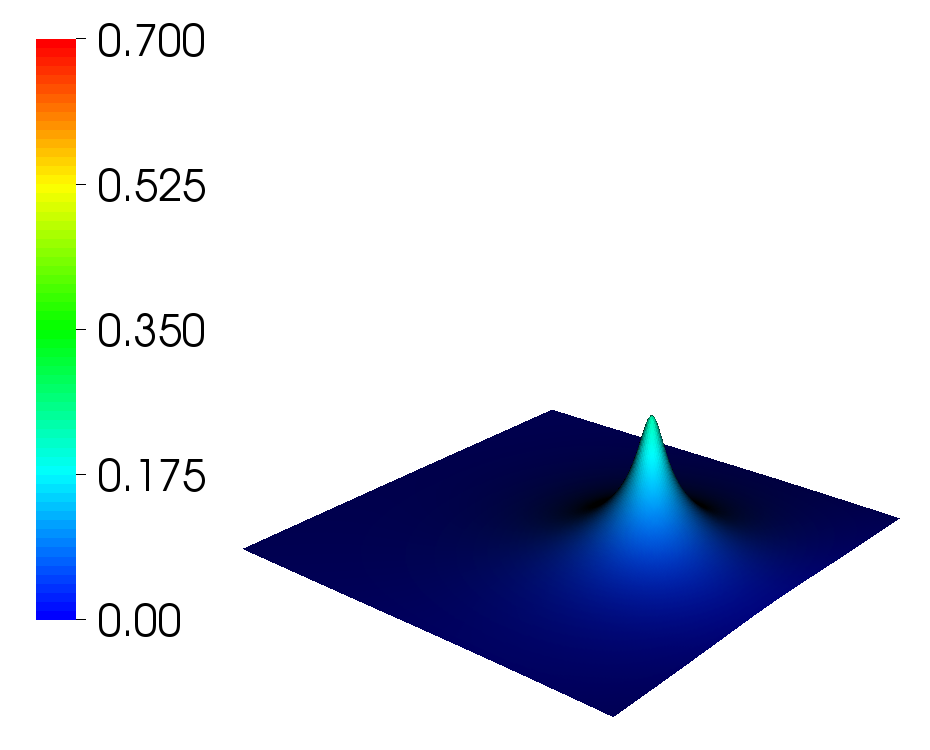}
  \end{minipage}
  \caption{Left: Incident light intensity $u$ for constant absorption
    coefficient. Center and right: Modulated light intensity $v^\xi$ for two
different focus
    points $\xi$. Note that $v$ depends on the focus position as well as the
    intensity of $u$ at the focus.}
  \label{fig:u-orig-v}
\end{figure}

Given $u$ and $|p^i|^2$, we compute the intensity of the modulated light
$v^i(x)$, using equations \eqref{eq:inverse_pr}. The equations are again
solved using
$\Q_1$ finite elements. Two examples for $v^i$ are shown in
Fig.~\ref{fig:u-orig-v} for two different focus positions.
The modulated light intensities $v^i$ are then evaluated at the sensor location
$\eta$ to yield the measurements $h(\xi_i) = v^i(\eta)$.

\subsection{Green's function and reconstruction}\label{sec:green-functions}

The reconstruction algorithm requires knowledge of the
Green's function $G$, which, given the absorption coefficient $\mu$ and
resulting diffusion coefficient $D$,
solves (\ref{eq:green}). Hence, we compute $G$ by solving another diffusion
problem with homogeneous Robin boundary conditions and a suitable approximation
to the delta function on the right hand side. As before, this is done using a
finite element scheme, where we choose a different, coarser mesh than in
forward problem calculations to avoid committing inverse crimes.

An obvious problem in the reconstruction formula (\ref{eq:mu-recon}) is that it
involves derivatives of the measurement data $h(\xi)$, which causes
instabilities in the presence of noise. Possible regularizations
for this problem are well-studied (e.g. \cite{Nam02}), and the stability
analysis in Section \ref{sec:stability} suggests that this is the only source
of instability in the reconstruction process. Hence, we opt not to add
extra regularization and compute the derivatives by a simple central
finite differencing scheme. Without adding noise to the measurements, it
turned out that in all of our computational experiments, the regularization
stemming from discretization on a fixed grid was sufficient for convergence of
the iterative scheme based on (\ref{eq:mu-recon}).

\subsection{Numerical phantoms}
\begin{figure}[htbp]
  \centering
  \begin{minipage}[b]{0.32\textwidth}
    \includegraphics[width=\textwidth]{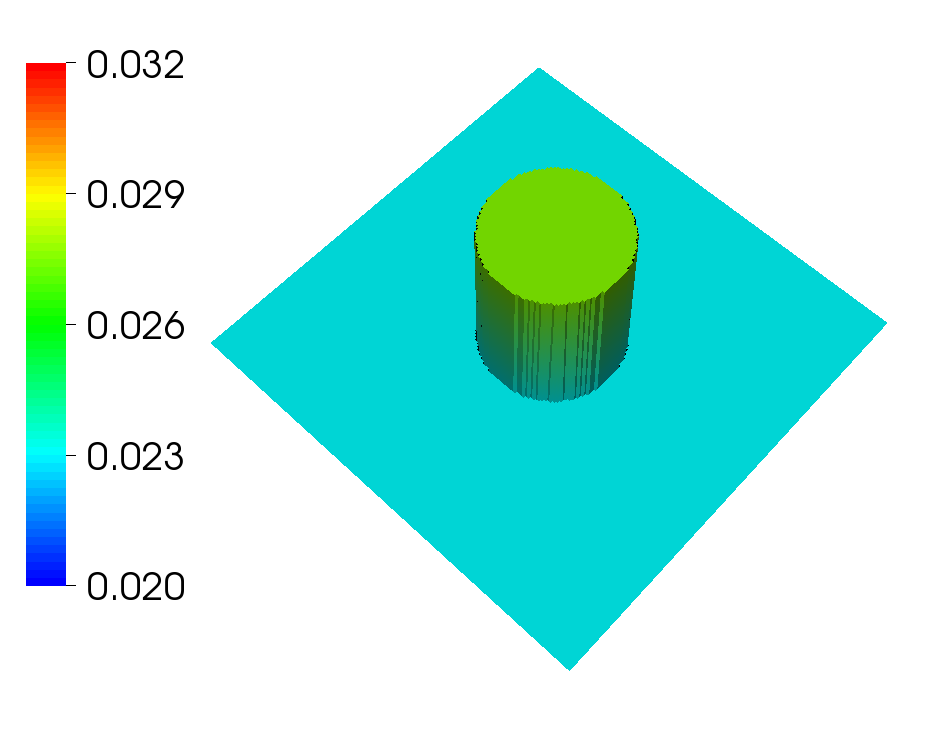}
  \end{minipage}
  \begin{minipage}[b]{0.32\textwidth}
    \includegraphics[width=\textwidth]{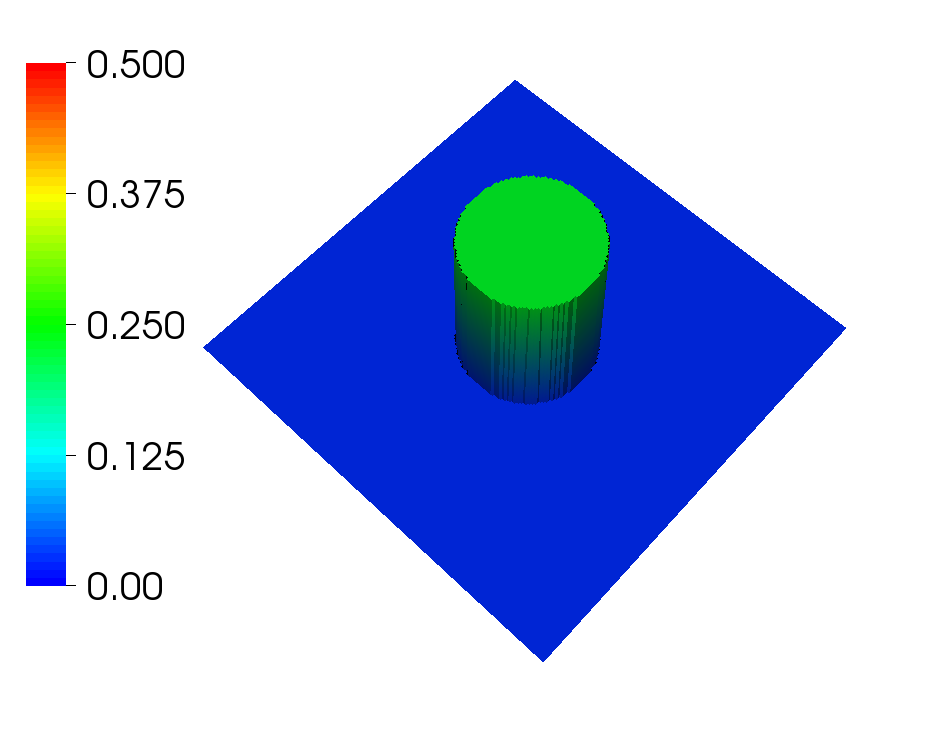}
  \end{minipage}
  \begin{minipage}[b]{0.32\textwidth}
  \includegraphics[width=\textwidth]{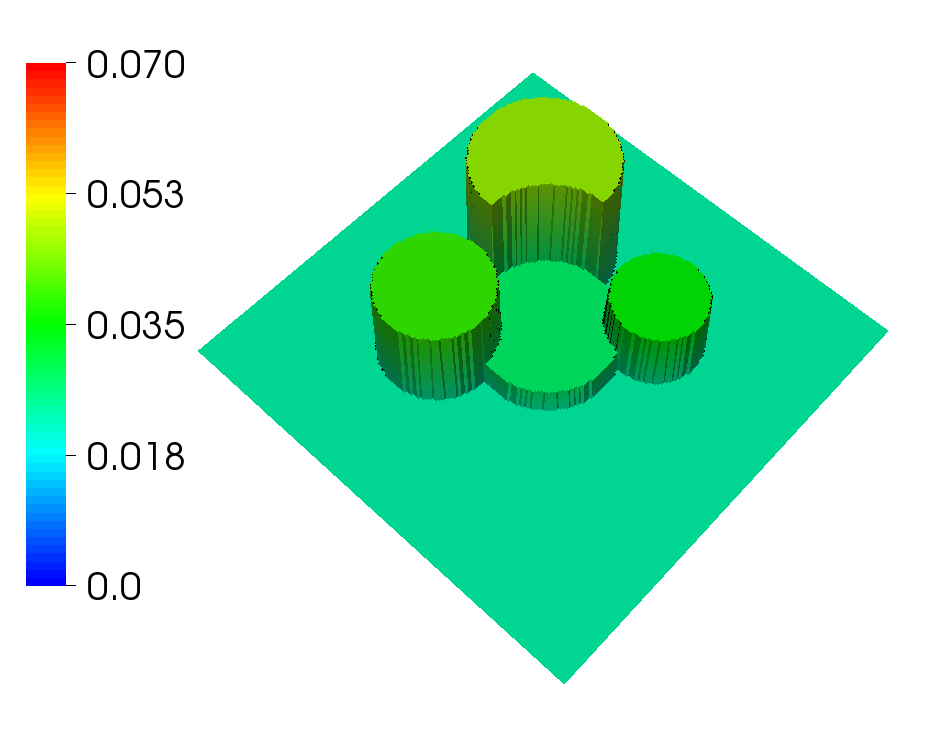}
  \end{minipage}
  \caption{Test cases for absorption coefficient $\mu^\ast$.}
  \label{fig:mu-star}
\end{figure}

To test our algorithms, we use three test cases in which the
true absorption coefficients have the following form:
\begin{itemize}
\item A disk-shaped inclusion $K\subset\Omega$ with midpoint $(2.5\rm{cm},
2.5\rm{cm})$ and radius $0.5\rm{cm}$. The absorption coefficient is assumed to
be equal to $\bar\mu$ outside the inclusion and slightly
higher inside:
\[
\mu^\ast (x)= \left\{
   \begin{array}[2]{l@{\;}l}
   \bar\mu &,\quad x \in \Omega\setminus K\\
   1.2\, \bar\mu &,\quad x \in K.
\end{array}\right.
\]

\item For the same inclusion $K$, a much higher absorption coefficient contrast
\[
 \mu^\ast (x)= \left\{
   \begin{array}[2]{l@{\;}l}
   \bar\mu &,\quad x \in \Omega\setminus K\\
   10\, \bar\mu &,\quad x \in K.
\end{array}\right.
\]

\item A more complicated coefficient with multiple inclusions of different
magnitude between $1.2\,\bar\mu$ and $2.0\,\bar\mu$. Their exact shape is
shown in Fig.~\ref{fig:mu-star}. This case tests the ability of our
algorithms to resolve several nearby objects.
\end{itemize}

For actual numerical values, we used $\bar\mu=0.023 \rm{cm^{-1}}$,
$\mu_s'=10.74 \rm{cm^{-1}}$ and $\gamma= 0.431 \rm{cm^{-1}}$ in our
computations. These values represent typical optical properties of soft tissue
\cite{MVD03}.

\subsection{Reconstruction results}
For the results shown in this section, measurements were produced
using the ultrasound signal arising from setting variances $\sigma_1 = \sigma_2
= 0.1\rm{cm}$ in the Gaussian (\ref{eq:field}), resulting in sharp focusing in
each direction (see the center panel of Fig.~\ref{fig:ultrasound} below).
Fig.~\ref{fig:reconstr} shows reconstructions of the three different
absorption coefficients for scanning the ultrasound focus $\xi^i$ on a
$100\!\times\!100$ mesh of points inside the area of interest $U$.

\begin{figure}[htbp]
  \centering
  \begin{minipage}[b]{0.32\textwidth}
    \includegraphics[width=\textwidth]{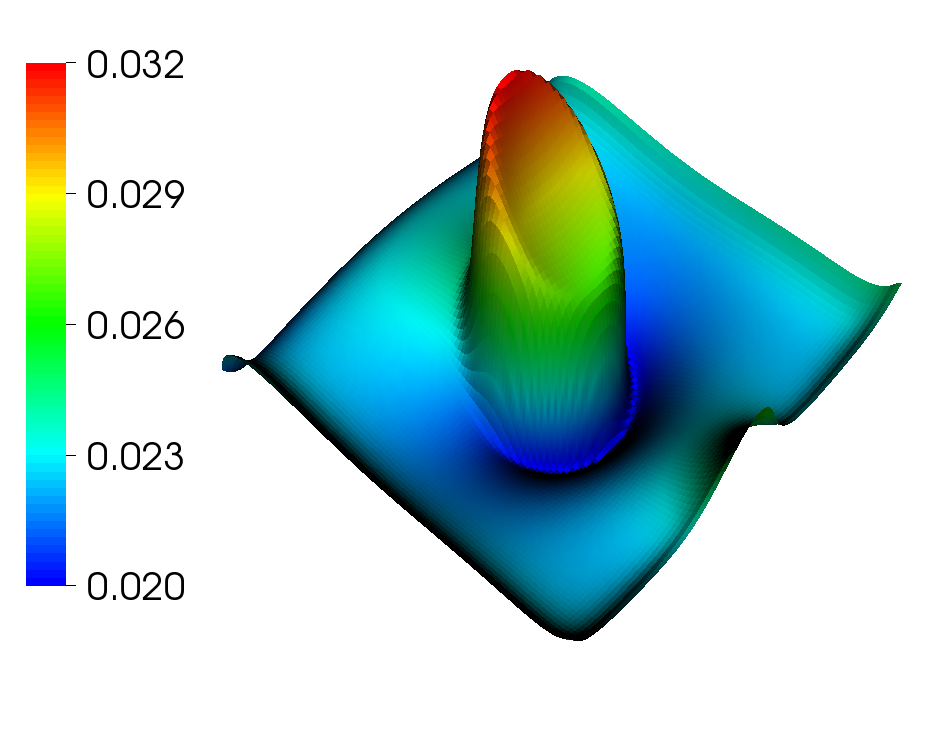}
  \end{minipage}
  \begin{minipage}[b]{0.32\textwidth}
    \includegraphics[width=\textwidth]{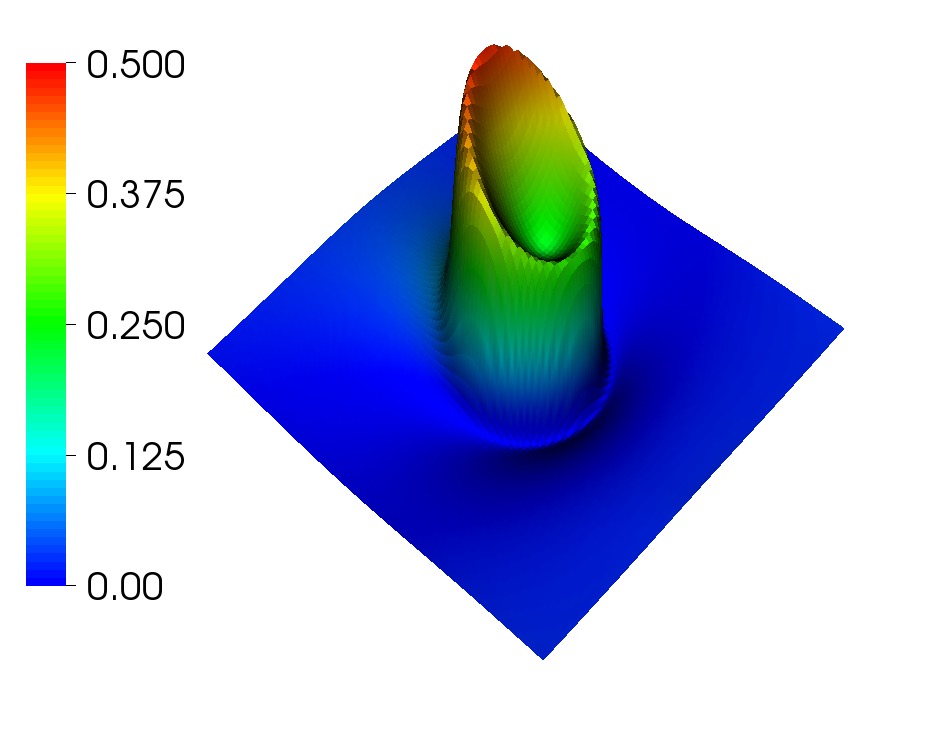}
  \end{minipage}
  \begin{minipage}[b]{0.32\textwidth}
    \includegraphics[width=\textwidth]{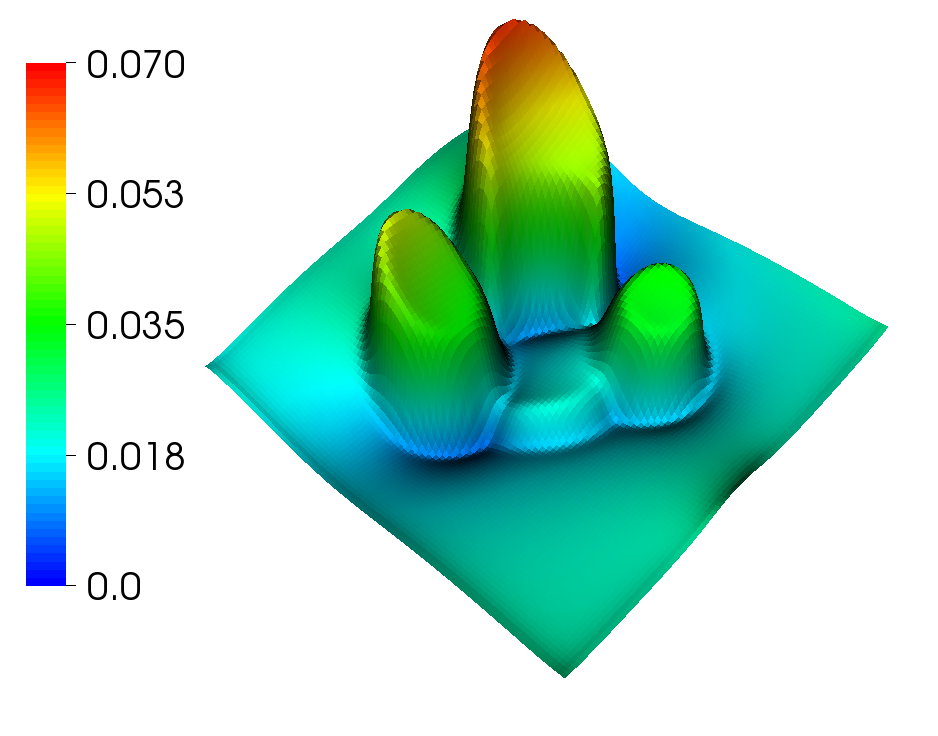}
  \end{minipage}
  \\
  \begin{minipage}[b]{0.32\textwidth}
  \includegraphics[width=\textwidth]{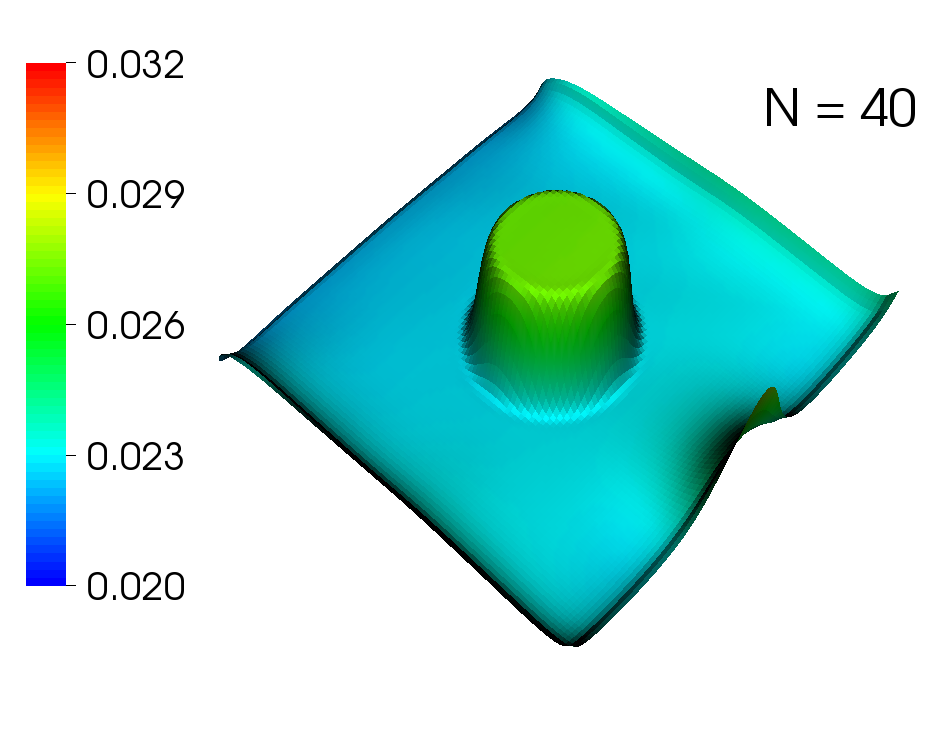}
  \end{minipage}
  \begin{minipage}[b]{0.32\textwidth}
  \includegraphics[width=\textwidth]{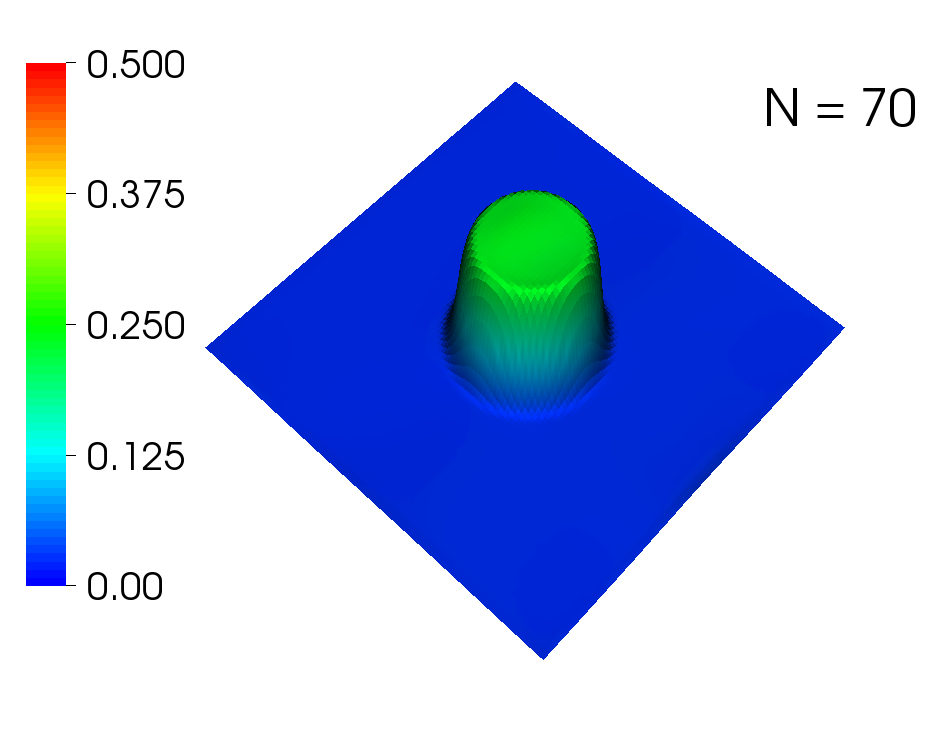}
  \end{minipage}
  \begin{minipage}[b]{0.32\textwidth}
  \includegraphics[width=\textwidth]{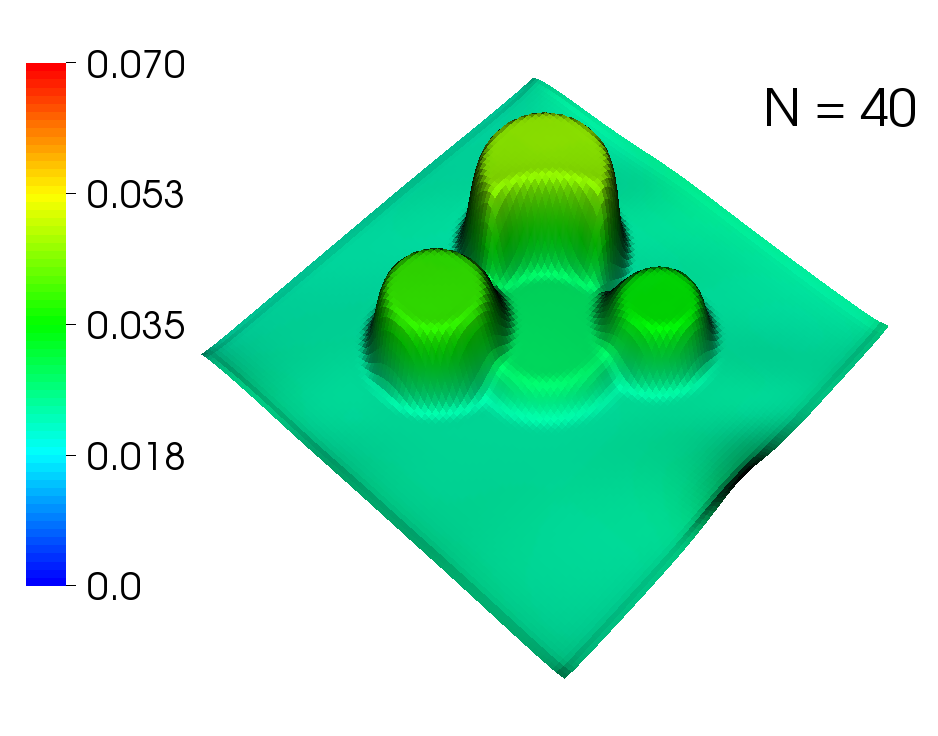}
  \end{minipage}
  \caption{Reconstruction results for the three coefficient cases: after the
first step of the algorithm (top) and after $N=40,70$ and $40$
iterations, respectively (bottom). }
  \label{fig:reconstr}
\end{figure}

The principal observation from these results is that under the main
assumptions of the model, i.e. turbid medium (and thus $\mu\ll \mu'_s$),
virtual light source, and strong focusing, our reconstruction scheme has four
desirable properties:
\begin{itemize}
\item It converges, even for the second case where (i) we start far away from
  the exact coefficient and (ii) the exact coefficient has a large dynamic
  range.
\item It is stable, i.e. the errors introduced through discretization of the
  equations, finite differencing of data, and using different meshes for
  reconstruction and generation of synthetic data do not lead to inaccurate
  reconstructions.
\item It can recover sharp interfaces without excessive blurring.
\item It can recover quantitatively correct values of absorption.
\end{itemize}
These are significant advantages compared to many other optical tomographic
methodologies.

\section{Stability of the linearized problem}
\label{sec:stability}

The quality of reconstructions shown above, especially the recovery of sharp
singularities, is
at first surprising, given that the standard OT problem is strongly ill-posed.
In this section, we will make a first step towards understanding the
stability of the UOT procedure.

Note that even though equations \eqref{eq:inverse_pr} defining
$u$ and $v$ are linear, the relation between the absorption coefficient
$\mu$ and the measurements $h$ is nonlinear. In
this section, we consider a (formal) linearization of the system
\eqref{eq:inverse_pr} that will allow us to gain some insight into the
local properties of the inverse problem.

Let $\Omega\subset\R^d$ with $d=2$ or $d=3$ be an open bounded domain with
$C^2$-boundary. We use a formal linearization, assuming that $\mu$ is a small
perturbation of a known absorption $\mu_0>0, \mu_0\in
C^{0,1}(\overline{\Omega})$, and then applying the formal asymptotic expansions
\begin{eqnarray*}
\mu(x) &\;=\;& \mu_0(x) + \eps \mu_1(x) + o(\eps),\\
 u(x) &\;=\;& u_0(x)+\eps u_1(x) + o(\eps), \\
v^\xi(x) &\;=\;& v_0^\xi(x) + \eps v_1^\xi(x) + o(\eps),
\end{eqnarray*}
where $\eps \rightarrow 0$. Our goal is to relate the
first order perturbations of the absorption coefficient $\mu_1$ and the
measurements $h_1(\xi):=v_1^\xi(\eta)$, where $\eta\in\partial\Omega$ is the
location of the detector.

Let us again assume perfectly focused ultrasound, i.e.
$|p^\xi(x)|^2=\delta(x-\xi)$. By inserting the above expansions into equations
\eqref{eq:inverse_pr} and sorting terms according to powers of $\eps$, we then
get the zeroth order perturbation system
\begin{eqnarray}
\label{eq:diffusion-u0}
  -\nabla\cdot D\nabla u_0(x) +\mu_0(x)
u_0(x)
  &=&0, 
  \\
\label{eq:diffusion-v0}
  -\nabla\cdot D\nabla v_0^\xi(x) +\mu_0(x)
v_0^\xi(x)
  &=&\alpha\delta (x-\xi) u_0(x), 
\end{eqnarray}
and the first order perturbation system
\begin{eqnarray}
\label{eq:diffusion-u1}
  -\nabla\cdot D\nabla u_1(x) +\mu_0(x)
  u_1(x)
  &=& - \mu_1(x) u_0(x), 
  \\
\label{eq:diffusion-v1}
  -\nabla\cdot D\nabla v_1^\xi(x) +\mu_0(x) v_1^\xi(x)
  &=&\alpha\delta (x-\xi) u_1(x) - \mu_1(x) v_0^\xi(x)
\end{eqnarray}
for all $x\in\Omega$, complemented by inhomogeneous Robin boundary
conditions as in \eqref{eq:diffusion-u-boundary} for $u_0$ and homogeneous Robin
boundary conditions for $v_0^\xi$, $u_1$ and $v_1^\xi$. Here we neglected the
(weak)
dependence of $D$ on $\mu$ and instead set $D\equiv \textrm{const}>0$ for the
rest of this section.

Equations \eqref{eq:diffusion-u0}--\eqref{eq:diffusion-v0} imply that
$u_0$ and $v_0^\xi$ are solutions to the forward model for absorption
coefficient $\mu_0$. The standard elliptic regularity theorems (e.g.,
\cite{GTr01}) imply $u_0\in H^3(\Omega)$, and by the
Sobolev embedding theorem $u_0\in C^{1}(\overline{\Omega})$ \cite{Evans98}.

Let us assume that the absorption coefficient is known near the boundary,
so
that it suffices to consider perturbations $\mu_1$ supported in an open set $U$
with $C^2$-boundary such that $\overline{U}\subset\Omega$. We assume
the data $h_1(\xi)$ to be given for all $\xi\in U$. In what follows, we
derive an explicit formula for the dependence of $\mu_1$ on $h_1$ and
then study properties of the corresponding linear operator.

Let us denote by $G_0(x, y)$ the Green's function as defined in \eqref{eq:green}
corresponding to the background absorption coefficient $\mu_0$. Equation
\eqref{eq:diffusion-v0} implies that for all $x \in\Omega$ and $\xi\in U$,
\begin{eqnarray*}
 v_0^\xi(x) \; &=& \; \int_\Omega \alpha G_0(x, z)
\delta(z - \xi) u_0(z) \; dz
\\
&=& \; \alpha G_0(x, \xi) u_0(\xi).
\end{eqnarray*}
From \eqref{eq:diffusion-v1} we can now deduce that
\begin{eqnarray*}
 v_1^\xi(x) \; &=& \; \int_\Omega G_0(x, z) \left[\alpha\delta(z - \xi)u_1(z) -
\mu_1(z) v_0^\xi(z)\right]\, dz
\\
&=& \; \alpha G_0(x, \xi) u_1(\xi) - \alpha u_0(\xi)\int_\Omega
G_0(x,z)G_0(z,\xi)
\mu_1(z)\; dz.
\end{eqnarray*}
Evaluating at $x=\eta$ and solving for $u_1$ yields
\[
u_1(\xi) \; = \; \frac{h_1(\xi)}{\alpha G_0(\eta,\xi)} +
\frac{u_0(\xi)}{G_0(\eta,\xi)} \int_\Omega G_0(\eta,z) G_0(z, \xi)
\mu_1(z)\; dz.
\]
We now use this expression to eliminate $u_1$ from \eqref{eq:diffusion-u1}.
Noting that the differential operators now act on $\xi$ and that
\[
\left[-\nabla_\xi\cdot D\nabla_\xi + \mu_0(\xi)\right]G_0(x, \xi) \; =
\; \delta(x- \xi),
\]
we get
\begin{align*}
  0
  &=
  u_0(\xi)\mu_1(\xi)
  +
  \left[-\nabla_\xi\cdot D\nabla_\xi + \mu_0(\xi)\right]
  \left(\frac{h_1(\xi)}{\alpha G_0(\eta,\xi)}\right)
  \\
  &\qquad+
  \left[-\nabla_\xi\cdot D\nabla_\xi + \mu_0(\xi)\right]
  \left(
    \frac{u_0(\xi)}{G_0(\eta,\xi)}
    \int_\Omega G_0(\eta,z) G_0(z,\xi) \mu_1(z) \; dz
  \right)
  \\
  &=
  u_0(\xi)\mu_1(\xi)
  +
  \left[-\nabla_\xi\cdot D\nabla_\xi + \mu_0(\xi)\right]
  \left(\frac{h_1(\xi)}{\alpha G_0(\eta,\xi)}\right)
  \\
  &\qquad+
  \left(
  \left[-\nabla_\xi\cdot D\nabla_\xi\right]
  \left[
    \frac{u_0(\xi)}{G_0(\eta,\xi)}
  \right]\right)
  \int_\Omega G_0(\eta,z) G_0(z,\xi) \mu_1(z) \; dz
  \\
  &\qquad-
  2 D\left[\nabla_\xi
    \left(
    \frac{u_0(\xi)}{G_0(\eta,\xi)}
    \right)\right]
  \cdot
  \left[\nabla_\xi
    \int_\Omega G_0(\eta,z) G_0(z,\xi) \mu_1(z) \; dz
    \right]
  \\
  &\qquad+
  \frac{u_0(\xi)}{G_0(\eta,\xi)}
  G_0(\eta,\xi) \mu_1(\xi).
\end{align*}

We will frequently view $G_0(\eta, y)$ as a function of $y$ in the following and
hence introduce the notation
\[
G_0^\eta(y) := G_0(\eta, y) \quad \textrm{for }y\in\overline{U}.
\]
Note that since $\eta\in\partial\Omega$, $G_0^\eta$ has no singularities on
$\overline{U}$ and hence is a regular solution to \eqref{eq:diffusion-u0}
there. The elliptic regularity and Sobolev embeddings
imply $G_0^\eta\in C^1(\overline{U})$.

Let us define the following operators acting on functions $g$ defined on $U$:
\begin{align}
\label{eq:K1}
K_1g(\xi) &:= -\frac{1}{2u_0(\xi)}\left(
  \left[-\nabla_\xi\cdot D\nabla_\xi\right]
  \left[
    \frac{u_0(\xi)}{G_0^\eta(\xi)}
  \right]\right)
  \int_U G_0^\eta(z) G_0(z,\xi) g(z) \; dz,
\\
\label{eq:K2}
K_2g(\xi) &:= \frac{D}{u_0(\xi)}\left[\nabla_\xi
    \left(
    \frac{u_0(\xi)}{G_0^\eta(\xi)}
    \right)\right]\cdot\left[\nabla_\xi
  \int_U G_0^\eta(z) G_0(z,\xi) g(z) \; dz\right],
\end{align}
and
$$
F := 1-K_1-K_2.
$$
In terms of these operators, our considerations above imply that $\mu_1$ is a
solution to the following linear equation:
\begin{equation}
\label{eq:fredholm}
F\mu_1(\xi)\;=\;-\frac{1}{2u_0(\xi)}\left[-\nabla_\xi\cdot
D\nabla_\xi +
\mu_0(\xi)\right]\left(\frac{h_1(\xi)}{\alpha G_0^\eta(\xi)}\right)
\end{equation}

In order for the above expressions to be well-defined, we have to make sure that
$u_0$ and $G_0^\eta$ are bounded away from zero on
$\overline{U}$. The following lemma follows
immediately from the Hopf Lemma (e.g., \cite{LiNi07,ProWein84}):

\begin{lemma}
\label{lem:u0}
There is a constant $c>0$ such that $u_0\ge c$ and  $G_0^\eta\ge c$ on
$\overline{U}$.
\end{lemma}

Next we consider the properties of the integral term involved in $K_1$ and
$K_2$. The important observation here is the following:
\begin{lemma}
The mapping
\begin{equation}
\label{eq:integral_op}
 g \mapsto \int_U G_0(z,\cdot) G_0^\eta(z) g(z)\, dz
\end{equation}
is a bounded linear operator from $L^2(U)$ to $H^2(U)$.
\end{lemma}
{\it Proof: } Let us assume that $g \in L^2(U)$. Since
$G_0^\eta\in C(\overline{U})$,  multiplication by $G_0^\eta$ is a bounded linear
operator on $L^2(U)$. The following integration against $G_0(z,\cdot)$ results
in the solution to the diffusion
equation with homogeneous Robin boundary condition and right hand side $G_0^\eta
g \in L^2(U)$. Elliptic regularity theory (e.g., \cite{GTr01,Evans98}) implies
that this is a continuous operator from $L^2(U)$ into $H^2(U)$.
\hfill$\scriptstyle{\Box}$

Because of the compact embedding of $H^2(U)$ in $L^2(U)$,
the operator defined by \eqref{eq:integral_op}, viewed as a mapping from
$L^2(U)$ to $L^2(U)$, is compact. In \eqref{eq:K1}, this operator is multiplied
by the factor
\begin{equation}
\label{eq:K1_factor}
-\frac{1}{2u_0(\xi)}\left(
  \left[-\nabla_\xi\cdot D\nabla_\xi\right]
  \left[
    \frac{u_0(\xi)}{G_0^\eta(\xi)}
  \right]\right).
\end{equation}
The functions $u_0, \nabla u_0, G_0^\eta$ and $\nabla G_0^\eta$ are all bounded
on $\overline{U}$ because $u_0, G_0^\eta\in C^1(\overline{U})$. Since $u_0$ and
$G_0^\eta$ satisfy \eqref{eq:diffusion-u0}, the terms $\nabla\cdot D \nabla u_0$
and $\nabla\cdot D \nabla G_0^\eta$ are bounded on $\overline{U}$ as well, and
$u_0^{-1}$ and $(G_0^\eta)^{-1}$ are bounded due to Lemma \ref{lem:u0}.
Consequently, multiplication by \eqref{eq:K1_factor} represents a bounded linear
operation on $L^2(U)$, and so $K_1$ is a compact operator in
$L^2(U)$. Similarly, $K_2$ is a compact operator in $L^2(U)$. This leads us to
the main result of this section:
\begin{theorem} $F:L^2(U) \to L^2(U)$ is a Fredholm operator of index zero.
\end{theorem}
Thus, the kernel ${\cal N}(F)$ of $F$ has finite dimension and the range
${\cal R}(F)$ is closed and of finite codimension, equal to the dimension of the
kernel. This immediately implies the following result:

\begin{corollary} $F$ as an operator from the quotient space $L^2(U) / {\cal
N}(F)$
to ${\cal R}(F)$ has bounded inverse, and the following norm equivalence
holds:
\begin{equation}\label{eq:stability}
 c_1 \|Ff\|_{L^2(U)} \;\le\;\|f\|_{L^2(U) / {\cal N}(F)} \;\le\; c_2
\|Ff\|_{L^2(U)}.
\end{equation}
\end{corollary}
The $L^2$-norm of the right hand side expression in \eqref{eq:fredholm} can be
estimated in terms of the $H^2$-norm of the measured perturbation $h_1$, so
that we obtain the following stability result:

\begin{theorem}\label{T:stability}
Under the stated assumptions, there is a constant $C>0$ such that the
following relation holds:
\begin{equation}
\label{eq:stability1}
 \|\mu_1\|_{L^2(U) / {\cal N}(F)} \;\le\;
C \|h_1\|_{H^2(U)}.
\end{equation}
\end{theorem}

We conjecture that the kernel ${\cal N}(F)$ is in fact trivial, and thus
the operator $F$ is invertible. This would imply that $\mu_1$ is
uniquely determined by the measured perturbation $h_1$, and allow us to replace
the quotient space norms in (\ref{eq:stability}) and (\ref{eq:stability1}) with
the regular $L^2$ norms. However, we have not been able to prove this result
yet.

Smoother norm coercive estimates for the absorption can be obtained if more is
assumed about the unperturbed absorption $\mu_0$ and the domain. For
instance, if $\mu_0\in C^\infty(\Omega)$, $S\in
C^\infty(\partial\Omega)$, and $\Omega$ has smooth boundary, the operators $K_1$
and $K_2$ defined in
\eqref{eq:K1}--\eqref{eq:K2}, are of order $-2$ and $-1$, respectively, in the
Sobolev scale:
\begin{eqnarray*}
K_1\;&:&\; H^s(U) \to H^{s+2}(U),\\
K_2\;&:&\; H^s(U) \to H^{s+1}(U).
\end{eqnarray*}
This and the Sobolev embedding theorem \cite{Ad75} imply that for any $s\ge 0$,
$F$ is Fredholm as an operator
\[
F\;:\; H^s(U)\to H^s(U).
\]
This, in turn, leads to the estimate
\begin{equation}
 \|f\|_{H^s(U)} \;\le\; c\left( \|f\|_{L^2(U)} +
\|Ff\|_{H^s(U)}\right)
\end{equation}
for all $f\in H^s(U)$. Thus, we have the following result:

\begin{theorem}\label{T:stability2}
Under the stated assumptions, for any $s>0$ there is a constant $C$ such
that
\[
\|\mu_1\|_{H^s(U)} \;\le\; C\left( \|\mu_1\|_{L^2(U)} +
\|h_1\|_{H^{s+2}(U)}\right).
\]
\end{theorem}
Clearly, if only a specific value of $s$ is of interest, the smoothness
assumptions on $\mu_0, S$ and $\partial\Omega$ can be relaxed appropriately.

\section{Conclusion and outlook}\label{sec:remarks}

In this paper, we have introduced a partial differential equation model of
ultrasound modulated optical tomography to derive a simple reconstruction
scheme for recovering the spatially varying absorption coefficient from boundary
measurements. While we could demonstrate stable, sharp and quantitatively
accurate reconstructions, some of the assumptions made here need to or can be
improved upon for practical applications. In particular, these are:

\paragraph{Detector locations}

In the discussion of stability above, as well as in our numerical
reconstructions, we have chosen a single detector point $\eta$. However, using
detectors distributed over a part $\Gamma$ of the boundary $\partial\Omega$
should help to suppress the effect of noise in the measured data.

\paragraph{Ultrasound signal with elongated focus}\label{sec:elongated}

\begin{figure}[htbp]
  \centering
  \begin{minipage}[b]{0.3\textwidth}
    \includegraphics[width=0.8\textwidth]{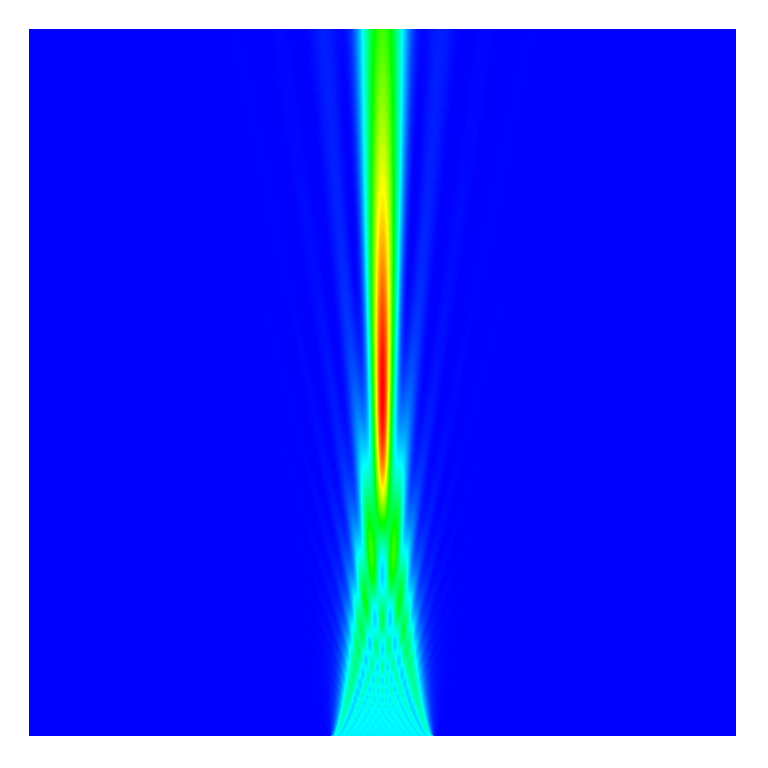}
  \end{minipage}
  \begin{minipage}[b]{0.3\textwidth}
    \includegraphics[width=0.8\textwidth]{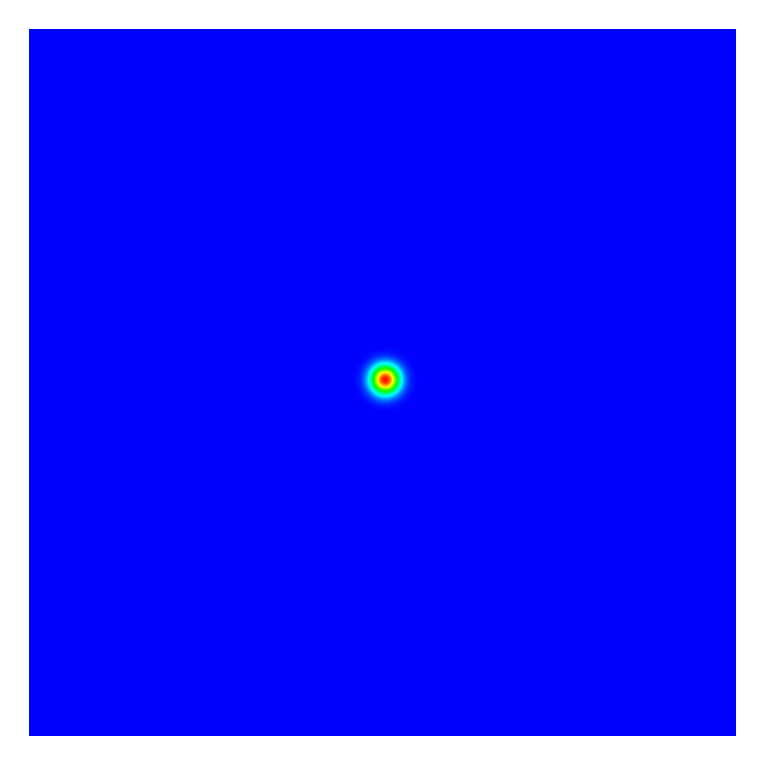}
  \end{minipage}
  \begin{minipage}[b]{0.3\textwidth}
    \includegraphics[width=0.8\textwidth]{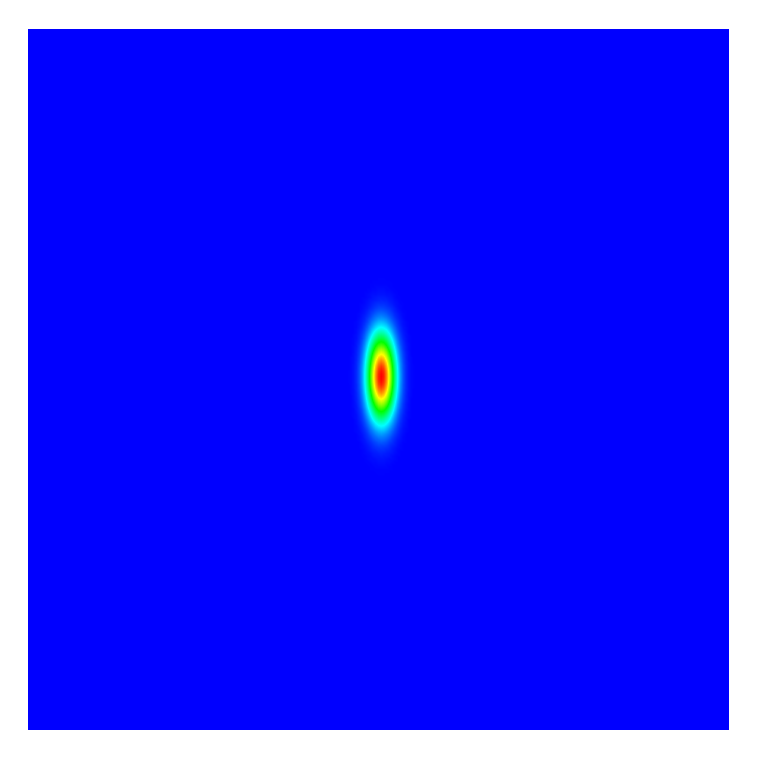}
  \end{minipage}
  \caption{Left: Simulated ultrasound pressure field $|p|^2$ with transducer at
the bottom. Middle: Gaussian ultrasound signal $|p|^2$ with $\sigma_1 =
\sigma_2 = 0.1$. Right: Gaussian signal with $\sigma_1 = 0.1$, $\sigma_2 = 0.3$.
}
  \label{fig:ultrasound}
\end{figure}

In practice, perfect focusing of ultrasound waves is not a
realistic assumption \cite{Localized}. How well an ultrasound
signal can be focused depends, in particular, on the geometry and bandwidth of
the transducer. For example, it is known from experimental measurements
(e.g., \cite{LevF01})
that focused ultrasound signals have an intensity profile similar to the one
shown in Fig.~\ref{fig:ultrasound} (left).
This signal has significantly sharper focus in the direction transverse to the
transducer lens, while the well-focused Gaussian signal used in our results
does not reflect this behavior.

To illustrate the effect of relaxing the assumption of perfect focus, we
computed reconstructions for the case where the ultrasound intensity is a
Gaussian signal with sharp focus in $x$-direction and elongated focus in
$y$-direction (Fig.~\ref{fig:ultrasound}, right).  As in the previous section,
the ultrasound focus $\xi^i$ is scanned on a $100\!\times\!100$ mesh to
produce synthetic measurements. At the same time, the reconstruction algorithm
is left unchanged, i.e. still assumes perfect focus.

\begin{figure}[tbp]
  \centering
  \begin{minipage}[b]{0.32\textwidth}
    \includegraphics[width=\textwidth]{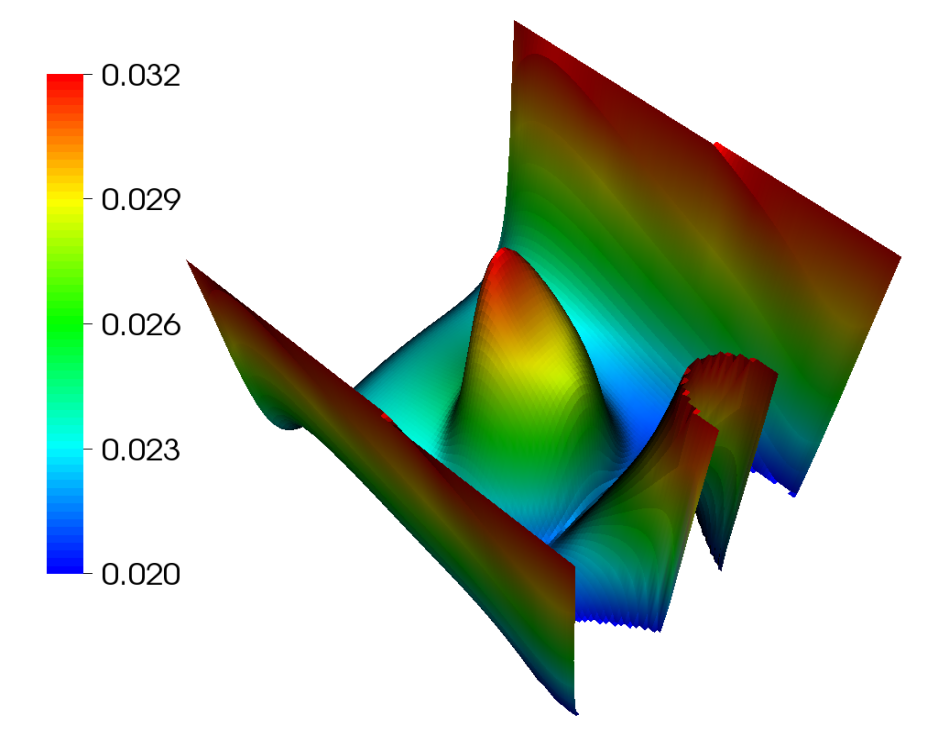}
  \end{minipage}
  \begin{minipage}[b]{0.32\textwidth}
    \includegraphics[width=\textwidth]{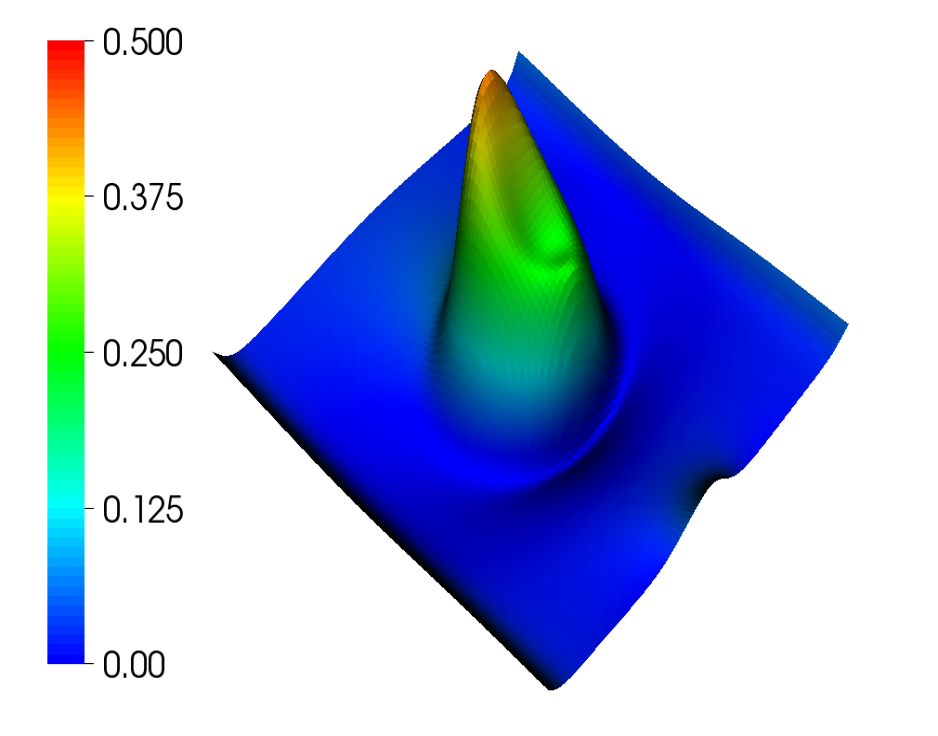}
  \end{minipage}
  \begin{minipage}[b]{0.32\textwidth}
    \includegraphics[width=\textwidth]{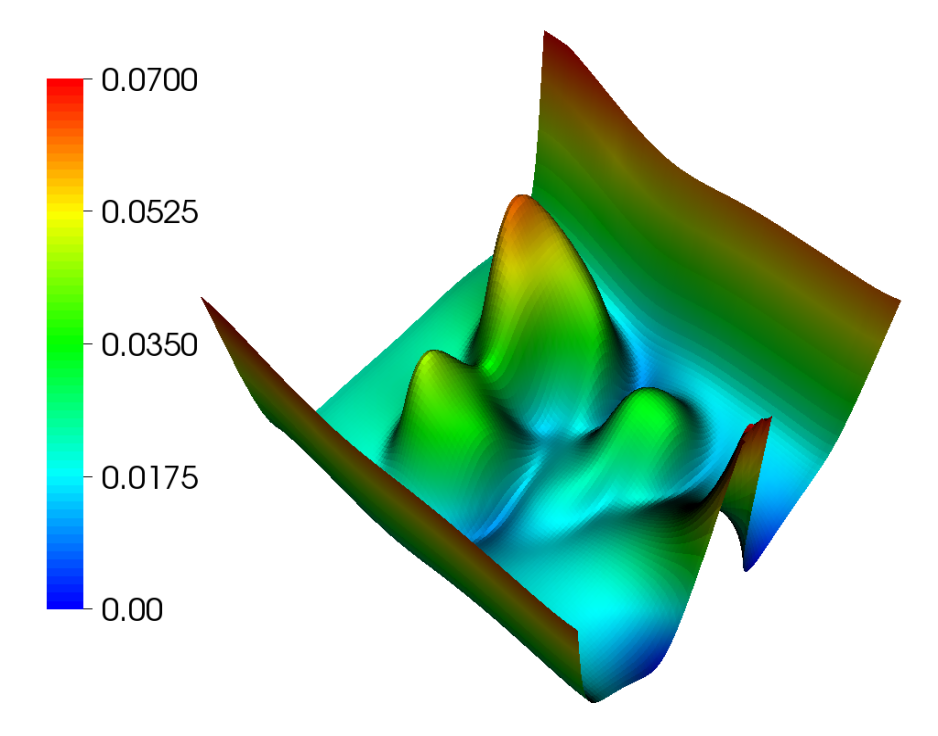}
  \end{minipage}
  \\
  \begin{minipage}[b]{0.32\textwidth}
  \includegraphics[width=\textwidth]{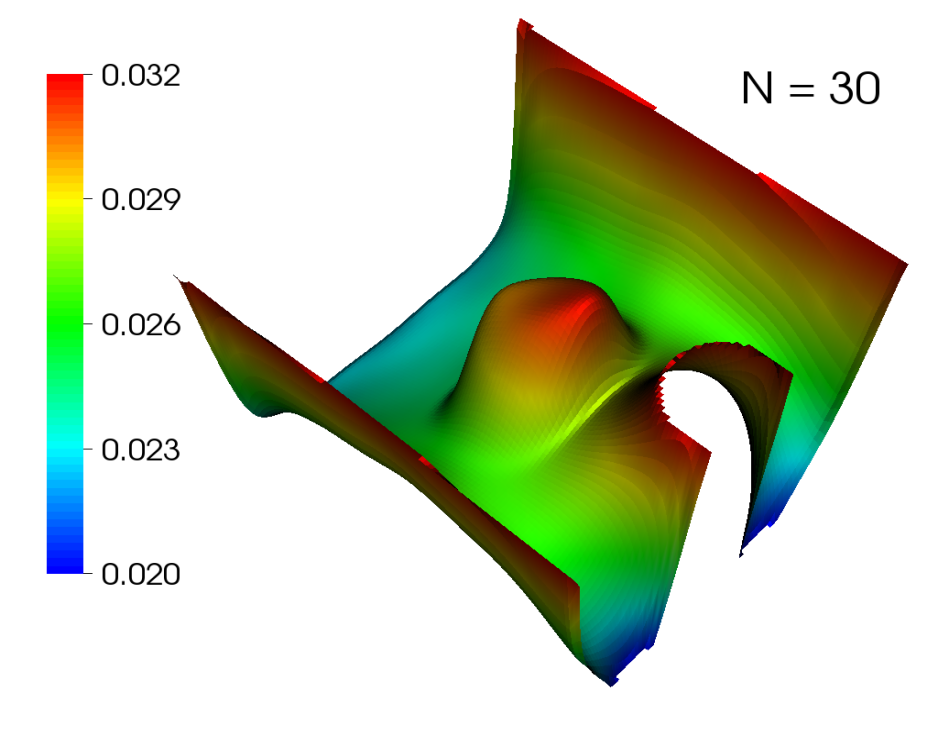}
  \end{minipage}
  \begin{minipage}[b]{0.32\textwidth}
  \includegraphics[width=\textwidth]{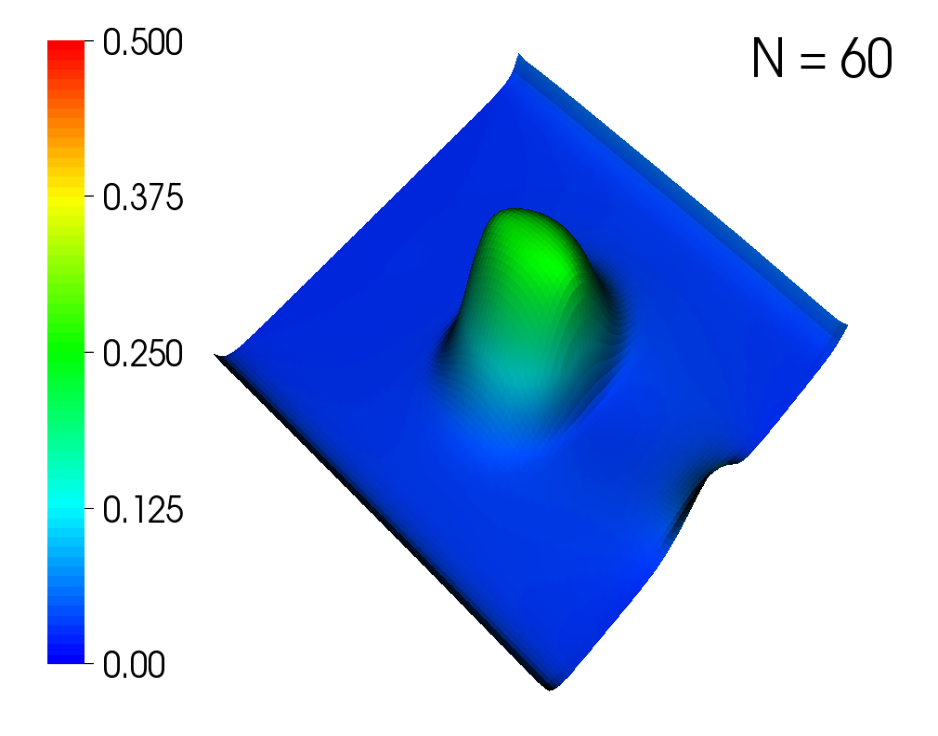}
  \end{minipage}
  \begin{minipage}[b]{0.32\textwidth}
  \includegraphics[width=\textwidth]{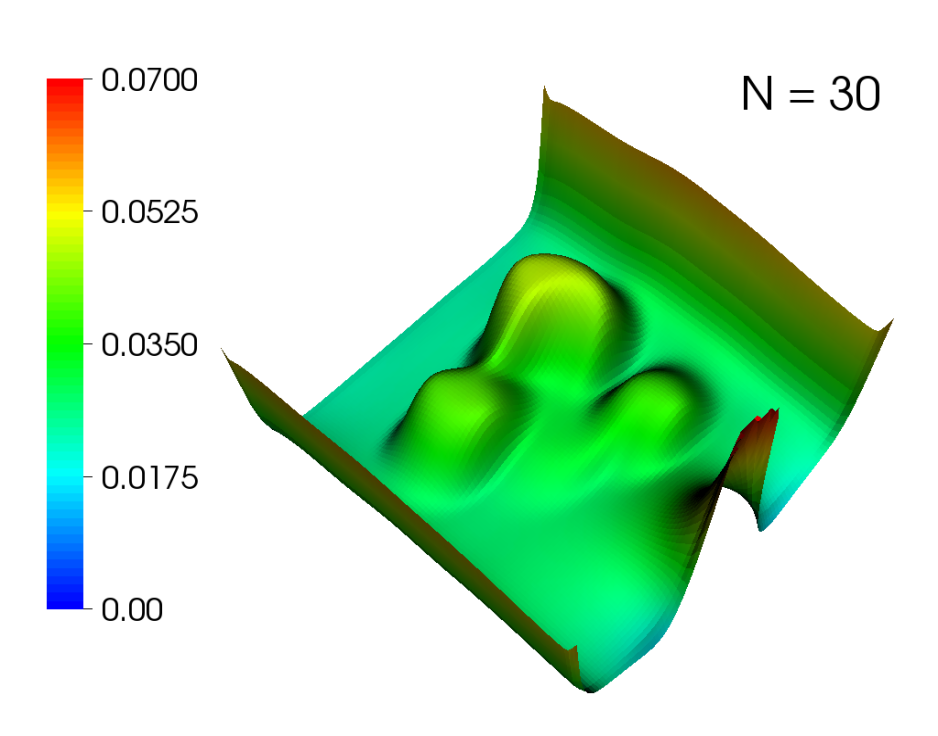}
  \end{minipage}
  \caption{Reconstruction results for ultrasound signal with elongated focus:
after the first step of the algorithm (top) and after $N$
iterations (bottom).}
  \label{fig:reconstr_sigma_y}
\end{figure}

Reconstruction results are shown in Fig.~\ref{fig:reconstr_sigma_y}. The
deterioration of the reconstruction -- in particular in the direction of the
ultrasound beam -- is obvious. The results also contain artifacts at the
vertical boundaries and close to the detector location. A more sophisticated
reconstruction scheme might be needed to treat the non-perfect focusing in these
calculations.

\paragraph{Synthetic focusing}\label{sec:synthetic}
Instead of attempting to perfectly focus the ultrasound waves in space, 
\emph{synthetic focusing} allows the use of non-localized
ultrasound fields and reconstructs the signal by superposition. This approach
was suggested in \cite{KuKu10}: It combines various basis sets of
non-focused ultrasound waves (e.g., spherical or monochromatic planar ones),
with a post-processing step that synthesizes the would-be response to a focused
illumination. In particular, in the case of spherical waves, the post-processing
(synthetic focusing procedure) is essentially equivalent to thermoacoustic
tomography inversion (see \cite{KuKu}). We plan to investigate the
applicability of this approach to UOT in the future.

\paragraph{Uniqueness of reconstruction}
Proving uniqueness of reconstruction, both in the non-linear and linearized
versions, still remains a challenge. In particular, we conjecture that the
operator $F$ in (\ref{eq:fredholm}) is in fact invertible, and thus there is
uniqueness of solution of the linearized problem, which would replace the
quotient space norms in (\ref{eq:stability}) and (\ref{eq:stability1}) with
the regular $L^2$ norms. At the same time, a complete characterization of the
kernel of the operator $F$ is non-trivial and left for future work.

\paragraph{Summary}
Despite these opportunities for future work,
in this paper, a diffusion based model is provided for the ultrasound
modulated optical tomography procedure using well focused ultrasound waves. An
iterative algorithm is suggested to recover absorption from measurements of the
amplitude of ultrasound modulation. The provided numerical results
show feasibility of the algorithm and possibility of
good reconstructions, both with regard to locating sharp interfaces, as well as
recovering correct numerical values of the absorption coefficient. Such
stability and resolution are impossible to achieve in standard optical
tomography. The stability of reconstructions is
explained by the stability estimates derived in Theorems \ref{T:stability} and
\ref{T:stability2} for a linearized model.

\section*{Acknowledgments}

The work of both authors was partially supported by NSF grant DMS-0604778 and
Award No.~KUS-C1-016-04 made by King Abdullah University of Science and
Technology (KAUST). The work of the second author was also partially supported
by U.S. Department of Energy grant DE-FG07-07ID14767 and by an Alfred P.~Sloan
Research Fellowship.
We wish to express our gratitude to these sources of support. We also
thank Prof.~P.~Kuchment, who suggested the approach for proving linear stability
in Section \ref{sec:stability}.

\bibliographystyle{model1-num-names}
\bibliography{umot}

\end{document}